\documentclass[12pt]{article}

\usepackage{amsfonts, amssymb, amsmath, amsthm}

\theoremstyle{plain}

\newtheorem{theorem}{Theorem}
\newtheorem{lemma}{Lemma}

\newtheorem{condition}{Condition}

\theoremstyle{definition}

\newtheorem{definition}{Definition}

\renewcommand{\Im}{{\rm Im\,}}

\newcommand{\supp}{{\rm supp\,}}

\newcommand{\dist}{{\rm dist}}

\begin{document}

%
%

\begin{center}
{\bf\Large Smoothness of Generalized Solutions for Nonlocal
Elliptic Problems on the Plane}

Pavel Gurevich\footnote{This research was supported by Russian
Foundation for Basic Research (grant No.~02-01-00312), Russian
Ministry for Education (grant No~E02-1.0-131), and INTAS
(grant~YSF~2002-008).}
\end{center}

\abstract{We study smoothness of generalized solutions of nonlocal
elliptic problems in plane bounded domains with piecewise smooth
boundary. The case where the support of nonlocal terms can
intersect the boundary is considered. We announce conditions that
are necessary and sufficient for any generalized solution to
possess an appropriate smoothness (in terms of Sobolev spaces).
The proofs are given in the forthcoming paper.}

\bigskip

{\bf 1.} The most difficult situation in the theory of elliptic
problems with nonlocal boundary-value conditions is that where the
support of nonlocal terms can intersect a boundary of a
domain~\cite{BitzDAN85}--\cite{SkRJMP}. In that case, solutions
can have power-law singularities near some points on the boundary.
In the present paper, we find out conditions that are necessary
and sufficient for any generalized solution $u\in W_2^1(G)$ of a
nonlocal problem in a plane bounded domain $G$ to belong to
$W_2^2(G)$. We study the case in which different nonlocal
conditions are set on different parts of the boundary,
coefficients of nonlocal terms supported near the points of
conjugation of boundary conditions are variable, and nonlocal
operators corresponding to nonlocal terms supported outside the
conjugation points are abstract. Both homogeneous and
nonhomogeneous nonlocal conditions are investigated. We consider a
nonlocal perturbation of the Dirichlet problem for an elliptic
equation of order two. However, the obtained results can be
generalized to elliptic equations of order $2m$ with general
nonlocal conditions.

\medskip

Let $G\subset{\mathbb R}^2$ be a bounded domain with boundary
$\partial G$. Introduce a set ${\mathcal K}\subset\partial G$
consisting of finitely many points. Let $\partial
G\setminus{\mathcal K}=\bigcup_{i=1}^{N}\Gamma_i$, where
$\Gamma_i$ are open (in the topology of $\partial G$)
$C^\infty$-curves. We assume that the domain $G$ is a plane angle
in some neighborhood of each point $g\in{\mathcal K}$. Denote by
${\bf P}$ a differential operator of order two, with smooth
complex-valued coefficients, properly elliptic in $\overline{G}$.

For any closed set $\mathcal M$, write $\mathcal
O_{\varepsilon}(\mathcal M)=\{y\in \mathbb R^2: \dist(y, \mathcal
M)<\varepsilon\}$, where ${\varepsilon}>0$.

We now define operators corresponding to nonlocal conditions near
the set $\mathcal K$. Let $\Omega_{is}$ ($i=1, \dots, N;$ $s=1,
\dots, S_i$) denote $C^\infty$-diffeomorphisms taking a
neighborhood ${\mathcal O}_i$ of the curve
$\overline{\Gamma_i\cap\mathcal O_{{\varepsilon}}(\mathcal K)}$
onto the set $\Omega_{is}({\mathcal O}_i)$ in such a way that
$\Omega_{is}(\Gamma_i\cap\mathcal O_{{\varepsilon}}(\mathcal
K))\subset G$ and $ \Omega_{is}(g)\in\mathcal K$ for $
g\in\overline{\Gamma_i}\cap\mathcal K. $ Thus, the transformations
$\Omega_{is}$ take the curves $\Gamma_i\cap\mathcal
O_{{\varepsilon}}(\mathcal K)$ strictly inside the domain $G$ and
their end points $\overline{\Gamma_i}\cap\mathcal K$ to the end
points.

Let us specify the structure of the transformations $\Omega_{is}$
near the set $\mathcal K$. Denote by $\Omega_{is}^{+1} $  the
transformation $\Omega_{is}:{\mathcal O}_i\to\Omega_{is}({\mathcal
O}_i)$ and by $\Omega_{is}^{-1}:\Omega_{is}({\mathcal
O}_i)\to{\mathcal O}_i$ the transformation inverse to
$\Omega_{is}.$ The set of the points
$\Omega_{i_qs_q}^{\pm1}(\dots\Omega_{i_1s_1}^{\pm1}(g))\in{\mathcal
K}$ ($1\le s_j\le S_{i_j},\ j=1, \dots, q$) is called an {\em
orbit} of the point $g\in{\mathcal K}$. In other words, the orbit
of $g\in{\mathcal K}$ is formed by the points that can be obtained
by consecutively applying the transformations $\Omega_{is}^{\pm1}$
to $g$. We assume for simplicity that the set $\mathcal
K=\{g_1,\dots,g_N\}$ consists of one orbit only.

Let $\varepsilon$ be small so that there exist neighborhoods
$\mathcal O_{\varepsilon_1}(g_j)$ of the points $g_j\in\mathcal K$
satisfying the following conditions: (I) $ \mathcal
O_{\varepsilon_1}(g_j)\supset\mathcal O_{\varepsilon}(g_j) $, (II)
the boundary $\partial G$ is an angle in the neighborhood
$\mathcal O_{\varepsilon_1}(g_j)$, (III) $\overline{\mathcal
O_{\varepsilon_1}(g_j)}\cap\overline{\mathcal
O_{\varepsilon_1}(g_k)}=\varnothing$ for any $g_j,g_k\in\mathcal
K$, $k\ne j$, (IV) if $g_j\in\overline{\Gamma_i}$ and
$\Omega_{is}(g_j)=g_k,$ then ${\mathcal
O}_{\varepsilon}(g_j)\subset\mathcal
 O_i$ and
 $\Omega_{is}\big({\mathcal
O}_{\varepsilon}(g_j)\big)\subset{\mathcal
O}_{\varepsilon_1}(g_k).$

For each point $g_j\in\overline{\Gamma_i}\cap\mathcal K$, fix a
transformation $y\mapsto y'(g_j)$ which is a composition of the
shift by the vector $-\overrightarrow{Og_j}$ and the rotation
through some angle so that the set ${\mathcal
O}_{\varepsilon_1}(g_j)$ is taken onto a neighborhood ${\mathcal
O}_{\varepsilon_1}(0)$ of the origin, whereas $G\cap{\mathcal
O}_{\varepsilon_1}(g_j)$ and $\Gamma_i\cap{\mathcal
O}_{\varepsilon_1}(g_j)$ are taken to the intersection of a plane
angle $K_j=\{y\in{\mathbb R}^2: r>0, |\omega|<\omega_j\}$ with
${\mathcal O}_{\varepsilon_1}(0)$ and to the intersection of the
side $\gamma_{j\sigma}=\{y\in\mathbb R^2: \omega=(-1)^\sigma
\omega_j\}$ ($\sigma=1$ or $\sigma=2$) of the angle $K_j$ with
${\mathcal O}_{\varepsilon_1}(0)$, respectively. Here $(\omega,r)$
are the polar coordinates, $0<\omega_j<\pi$.

\begin{condition}\label{condK1}
The above change of variables $y\mapsto y'(g_j)$ for
$y\in{\mathcal O}_{\varepsilon}(g_j)$,
$g_j\in\overline{\Gamma_i}\cap\mathcal K$, reduces the
transformation $\Omega_{is}(y)$ to the composition of a rotation
and a homothety in the new variables $y'$.
\end{condition}

Introduce the nonlocal operators $\mathbf B_{i}^1$ by the formula
$
 \mathbf B_{i}^1u=\sum\limits_{s=1}^{S_i}
   b_{is}(y)u(\Omega_{is}(y))$,
$y\in\Gamma_i\cap\mathcal O_{\varepsilon}(\mathcal K)$, $\mathbf
B_{i}^1u=0$, $y\in\Gamma_i\setminus(\Gamma_i\cap\mathcal
O_{\varepsilon}(\mathcal K))$, where $b_{is}\in C^\infty(\mathbb
R^2)$ and  $\supp b_{is}\subset\mathcal O_{{\varepsilon}}(\mathcal
K)$. Since $\mathbf B_{i}^1u=0$ whenever $\supp u\subset\overline{
G}\setminus\overline{\mathcal O_{{\varepsilon_1}}(\mathcal K)}$,
we say that the operators $\mathbf B_{i}^1$ \textit{correspond to
nonlocal terms supported near the set} $\mathcal K$.

\smallskip

Consider the operators $\mathbf B_{i}^2$ satisfying the following
condition (cf.~(2.5), (2.6) in~\cite{SkMs86} and~(3.4), (3.5)
in~\cite{SkJMAA}).

\begin{condition}\label{condSeparK23}
There exist numbers $\varkappa_1>\varkappa_2>0$ and $\rho>0$ such
that the inequalities
\begin{equation}\label{eqSeparK23'}
  \|\mathbf B^2_{i}u\|_{W_2^{3/2}(\Gamma_i)}\le c_1
  \|u\|_{W_2^{2}(G\setminus\overline{\mathcal O_{\varkappa_1}(\mathcal
  K)})},\quad
  \|\mathbf B^2_{i}u\|_{W_2^{3/2}
   (\Gamma_i\setminus\overline{\mathcal O_{\varkappa_2}(\mathcal K)})}\le
  c_2 \|u\|_{W_2^{2}(G_\rho)}
\end{equation}
hold for any $u\in W_2^{2}(G\setminus\overline{\mathcal
O_{\varkappa_1}(\mathcal
  K)})\cap W_2^{2}(G_\rho)$, where $G_\rho=\{y\in G: \dist(y,
\partial G)>\rho\}$, $i=1, \dots, N$, $c_1,c_2>0$.
\end{condition}

In particular, the first inequality in~\eqref{eqSeparK23'} means
that $\mathbf B_{i}^2u=0$ whenever $\supp u\subset \mathcal
O_{\varkappa_1}(\mathcal K)$. Therefore, we say that the operators
$\mathbf B_{i}^2$ \textit{correspond to nonlocal terms supported
outside the set} $\mathcal K$. Examples of the operators $\mathbf
B_i^2$ can be found in~\cite{SkMs86,SkDu91}.

\smallskip

We assume that Conditions~\ref{condK1} and~\ref{condSeparK23} are
fulfilled throughout the paper.

\smallskip

Consider the following nonlocal elliptic boundary-value problem:
\begin{align}
 {\bf P} u=f_0(y) \quad &(y\in G),\label{eqPinG}\\
     u|_{\Gamma_i}+\mathbf B_{i}^1 u+\mathbf B_{i}^2 u=
   f_{i}(y)\quad
    &(y\in \Gamma_i;\ i=1, \dots, N).\label{eqBinG}
\end{align}

Denote $\mathcal W_2^{k-1/2}(\partial G)=\prod_{i=1}^{N}
W_2^{k-1/2}(\Gamma_i)$ for $k\in\mathbb N$. For any set $X\in
\mathbb R^2$ having a nonempty interior, we denote by
$C_0^\infty(X)$ the set of functions infinitely differentiable in
$\overline{ X}$ and supported in $X$.

\begin{definition}\label{defGenSol2}
A function $u\in W_2^{1}(G)$ is called a {\em generalized
solution} of problem~\eqref{eqPinG}, \eqref{eqBinG} with
right-hand side $\{f_0,f_i\}\in L_2(G)\times \mathcal
W_2^{1/2}(\partial G)$ if $u$ satisfies nonlocal
conditions~\eqref{eqBinG} (where the equalities are understood as
those in $W_2^{1/2}(\Gamma_i)$) and Eq.~\eqref{eqPinG} in the
sense of distributions.
\end{definition}

We now write a model nonlocal problem corresponding to the points
of the set (orbit) ${\mathcal K}$. Denote by $u_j(y)$ the function
$u(y)$ for $y\in{\mathcal O}_{\varepsilon_1}(g_j)$. If
$g_j\in\overline{\Gamma_i},$ $y\in{\mathcal
O}_{\varepsilon}(g_j),$  and $\Omega_{is}(y)\in{\mathcal
O}_{\varepsilon_1}(g_k),$ then denote by $u_k(\Omega_{is}(y))$ the
function $u(\Omega_{is}(y))$. In that case, nonlocal
problem~(\ref{eqPinG}), (\ref{eqBinG}) acquires the following form
in the $\varepsilon$-neighborhood of the set (orbit) $\mathcal K$:
$$
 {\bf P} u_j=f_0(y) \quad (y\in\mathcal O_\varepsilon(g_j)\cap
 G),
$$
$$
\begin{aligned}
u_j(y)|_{\mathcal O_\varepsilon(g_j)\cap\Gamma_i}+
\sum\limits_{s=1}^{S_i} b_{is}(y) u_k(\Omega_{is}(y))|_{\mathcal
O_\varepsilon(g_j)\cap\Gamma_i}
=\psi_i(y) \\
\big(y\in \mathcal O_\varepsilon(g_j)\cap\Gamma_i;\ i\in\{1\le
i\le N: g_j\in\overline{\Gamma_i}\};\ j=1, \dots, N\big),
\end{aligned}
$$
where $\psi_i=f_i-\mathbf B_{i}^2u$. Let $y\mapsto y'(g_j)$ be the
above change of variables. Set $K_j^\varepsilon=K_j\cap\mathcal
O_\varepsilon(0)$,
$\gamma_{j\sigma}^\varepsilon=\gamma_{j\sigma}\cap\mathcal
O_\varepsilon(0)$ and introduce the functions
\begin{equation}\label{eqytoy'}
U_j(y')=u_j(y(y')),\ F_j(y')=f_0(y(y')),\ y'\in
K_j^\varepsilon,\qquad \Psi_{j\sigma}(y')=\psi_i(y(y')),\
y'\in\gamma_{j\sigma}^\varepsilon,
\end{equation}
where $\sigma=1$ $(\sigma=2)$ if the transformation $y\mapsto
y'(g_j)$ takes $\Gamma_i$ to the side $\gamma_{j1}$
($\gamma_{j2}$) of the angle $K_j$. In what follows, we write $y$
instead of $y'$. Using Condition~\ref{condK1}, we can write
problem~(\ref{eqPinG}), (\ref{eqBinG}) as follows:
\begin{gather}
  {\bf P}_{j}U_j=F_{j}(y) \quad (y\in
  K_j^\varepsilon),\label{eqPinK}\\
  {\mathbf B}_{j\sigma}U\equiv\sum\limits_{k,s}
      b_{j\sigma ks}(y)U_k({\mathcal G}_{j\sigma ks}y)
    =\Psi_{j\sigma}(y) \quad (y\in\gamma_{j\sigma}^\varepsilon).\label{eqBinK}
\end{gather}
Here (and below unless otherwise stated) $j, k=1, \dots, N;$
$\sigma=1, 2;$ $s=0, \dots, S_{j\sigma k}$; $\mathbf P_j$ is an
elliptic differential operator of order two with smooth
coefficients; $U=(U_1,\dots,U_N)$; $b_{j\sigma ks}(y)$ are smooth
functions, $b_{j\sigma j0}(y)\equiv 1$; ${\mathcal G}_{j\sigma
ks}$ is the operator of rotation through an angle~$\omega_{j\sigma
ks}$ and of homothety with a coefficient~$\chi_{j\sigma ks}>0$ in
the~$y$-plane. Moreover, $|(-1)^\sigma \omega_{j}+\omega_{j\sigma
ks}|<\omega_{k}$ for $(k,s)\ne(j,0)$ and $\omega_{j\sigma j0}=0$,
$\chi_{j\sigma j0}=1$ (i.e., $\mathcal G_{j\sigma j0}y\equiv y$).

Write the principal parts of the operators $\mathbf P_j$ at the
point $y=0$ in polar coordinates, $r^{-2}\tilde{\mathcal
P_j}(\omega,\partial/\partial\omega,r\partial/\partial r)$.
Consider the analytic operator-valued function $\tilde{\mathcal
L}(\lambda):\prod_j W_2^{2}(-\omega_j,\omega_j)\to\prod_j
(L_2(-\omega_j, \omega_j)\times\mathbb C^2)$ given by $
\tilde{\mathcal L}(\lambda)\varphi=\big\{\tilde{\mathcal
P}_j(\omega,
\partial/\partial\omega, i\lambda)\varphi_j,\
  \sum\limits_{k,s} (\chi_{j\sigma ks})^{i\lambda}
              b_{j\sigma ks}(0)\varphi_k((-1)^\sigma \omega_j+\omega_{j\sigma ks})\big\}.
$ Main definitions and facts concerning analytic operator-valued
functions can be found in~\cite{GS}. It is fundamental that the
{\em spectrum of the operator $\tilde{\mathcal L}(\lambda)$ is
discrete and, for any numbers $c_1<c_2$, the band
$c_1<\Im\lambda<c_2$ contains at most finitely many eigenvalues of
the operator $\tilde{\mathcal L}(\lambda)$} (see~\cite{SkDu90}).
Spectral properties of the operator $\tilde{\mathcal L}(\lambda)$
play a crucial role in the study of smoothness of generalized
solutions.

\medskip

{\bf 2.} Let $\lambda=\lambda_0$ be an eigenvalue of the operator
$\tilde{\mathcal L}(\lambda)$.

\begin{definition}\label{defRegEigVal}
We say that $\lambda_0$ is a {\em proper eigenvalue} if none of
the eigenvectors $\varphi(\omega)=(\varphi_{1}(\omega),\dots,
\varphi_{N}(\omega))$ corresponding to $\lambda_0$ has an
associated vector, whereas the functions
$r^{i\lambda_0}\varphi_{j}(\omega)$, $j=1, \dots, N$, are
polynomials in $y_1, y_2$. An eigenvalue which is not proper is
said to be {\em improper}.
\end{definition}
The notion of proper eigenvalue was originally proposed by
Kondrat'ev~\cite{KondrTMMO67} for ``local'' boundary-value
problems in nonsmooth domains.

\begin{theorem}\label{thuinW_2^2NoEigen}
1. Let the band $-1\le\Im\lambda<0$ contain no eigenvalues of the
operator $\tilde{\mathcal L}(\lambda)$, and let $u\in W_2^1(G)$ be
a generalized solution of problem~\eqref{eqPinG}, \eqref{eqBinG}
with right-hand side $\{f_0,f_i\}\in L_2(G)\times \mathcal
W_2^{3/2}(\partial G)$. Then $u\in W_2^2(G)$.

2. Let the band $-1\le\Im\lambda<0$ contain an improper eigenvalue
of the operator $\tilde{\mathcal L}(\lambda)$. Then there exists a
generalized solution $u\in W_2^1(G)$ of problem~\eqref{eqPinG},
\eqref{eqBinG} with certain right-hand side $\{f_0,0\}$, $f_0\in
L_2(G)$, such that $u\notin W_2^2(G)$.
\end{theorem}

\medskip

It remains to study the case in which the following condition
holds.
\begin{condition}\label{condProperEigen}
The band $-1\le\Im\lambda<0$ contains a unique eigenvalue
$\lambda=-i$ of the operator $\tilde{\mathcal L}(\lambda)$, and
this eigenvalue is a proper one.
\end{condition}

We first consider problem~\eqref{eqPinG}, \eqref{eqBinG} with
nonhomogeneous nonlocal conditions.

Denote by $\tau_{j\sigma}$ the unit vector co-directed with the
ray~$\gamma_{j\sigma}$. Consider the operators $
  \dfrac{\partial} {\partial
  \tau_{j\sigma}} \Big(\sum\limits_{k,s}b_{j\sigma
  ks}(0)U_k({\mathcal G}_{j\sigma ks}y)\Big).
$ Using the chain rule, we write them as follows:
\begin{equation}\label{eqDiffB}
 \sum\limits_{k,s}(\hat B_{j\sigma ks}(D_y)U_k)({\mathcal G}_{j\sigma ks}y),
\end{equation}
where $\hat B_{j\sigma ks}(D_y)$ are first-order differential
operators with constant coefficients. In particular, we have $\hat
B_{j\sigma j0}(D_y)=
 {\partial}/
 {\partial \tau_{j\sigma}}$ because ${\mathcal G}_{j\sigma j0}y\equiv
 y$. Formally replacing the nonlocal operators by the corresponding local ones
in~(\ref{eqDiffB}), we introduce the operators
\begin{equation}\label{eqSystemB}
 \hat{\mathcal B}_{j\sigma}(D_y)U\equiv
 \sum\limits_{k,s}\hat B_{j\sigma ks}(D_y)U_k(y).
\end{equation}
If Condition~\ref{condProperEigen} holds, then the system of
operators~\eqref{eqSystemB} is linearly
dependent~\cite{GurRJMP03}. Let
\begin{equation}\label{eqSystemB'}
\{\hat{\mathcal B}_{j'\sigma'}(D_y)\}
\end{equation}
be a maximal linearly independent subsystem of
system~\eqref{eqSystemB}. In that case, any operator
$\hat{\mathcal B}_{j\sigma}(D_y)$ which does not enter
system~\eqref{eqSystemB'} can be represented as follows:
\begin{equation}\label{eqBviaB'}
\hat{\mathcal
B}_{j\sigma}(D_y)=\sum\limits_{j',\sigma'}\beta_{j\sigma}^{j'\sigma'}\hat{\mathcal
B}_{j'\sigma'}(D_y),
\end{equation}
where $\beta_{j\sigma}^{j'\sigma'}$ are some constants. Let
$Z_{j\sigma}\in W_2^{3/2}(\gamma_{j\sigma}^\varepsilon)$ be
arbitrary functions. Set $
Z^0_{j\sigma}(r)=Z_{j\sigma}(y)|_{y=(r\cos\omega_j,\,
r(-1)^\sigma\sin\omega_j)}$. It is clear that $Z^0_{j\sigma}\in
W_2^{3/2}(0,\varepsilon)$.

\begin{definition}
Let $\beta_{j'\sigma'}$ be the constants occurring
in~\eqref{eqBviaB'}. If the relations
\begin{equation}\label{eqConsistencyZ}
\int\limits_{0}^\varepsilon
r^{-1}\Bigg|\frac{d}{dr}\bigg(Z^0_{j\sigma}-\sum\limits_{j',\sigma'}\beta_{j\sigma}^{j'\sigma'}Z^0_{j'\sigma'}\bigg)\Bigg|^2dr<\infty
\end{equation}
hold for all indices $j,\sigma$ corresponding to the operators of
system~\eqref{eqSystemB} which do not enter
system~\eqref{eqSystemB'}, then we say that the {\em functions
$Z_{j\sigma}$ satisfy the consistency
condition~\eqref{eqConsistencyZ}}.
\end{definition}

Let us formulate conditions which ensure that generalized
solutions are smooth. We first show that right-hand sides $f_i$ in
nonlocal conditions~\eqref{eqBinG} cannot be arbitrary functions
from the space $W_2^{3/2}(\Gamma_i)$.

Consider the change of variables $y\mapsto y'(g_j)$ described in
Sec.~1. Introduce the functions
$$
F_{j\sigma}(y')=f_i(y(y')),\quad y'\in\gamma_{j\sigma}^\varepsilon
$$
(cf. functions~\eqref{eqytoy'}). Denote by $\mathcal
S_2^{3/2}(\partial G)$ the set consisting of functions
$\{f_i\}\in\mathcal W_2^{3/2}(\partial G)$ such that the functions
$F_{j\sigma}$ satisfy the consistency
condition~\eqref{eqConsistencyZ}. The set $\mathcal
S_2^{3/2}(\partial G)$ is not closed in the topology of $\mathcal
W_2^{3/2}(\partial G)$ (see.~\cite[Lemma~3.2]{GurRJMP03}).

\begin{lemma}\label{lUNonSmFNonConsist}
Let Condition~$\ref{condProperEigen}$ hold. Then there exist a
function $\{f_0,f_i\}\in L_2(G)\times \mathcal W_2^{3/2}(\partial
G)$, $\{f_i\}\notin\mathcal S_2^{3/2}(\partial G)$, and a function
$u\in W_2^1(G)$ such that $u$ is a generalized solution of
problem~\eqref{eqPinG}, \eqref{eqBinG} with right-hand side
$\{f_0,f_i\}$ and $u\notin W_2^2(G)$.
\end{lemma}

It follows from Lemma~\ref{lUNonSmFNonConsist} that, if one wants
{\em any} generalized solution of problem~\eqref{eqPinG},
\eqref{eqBinG} be smooth, then one must take right-hand sides
$\{f_0,f_i\}$ from the space $L_2(G)\times \mathcal
S_2^{3/2}(\partial G)$.

\medskip

Let $v\in W_2^{2}(G\setminus\overline{\mathcal
O_{\varkappa_1}(\mathcal
  K)})$ be an arbitrary function.
Consider the change of variables $y\mapsto y'(g_j)$ from Sec.~1
again and introduce the functions
$$
B^v_{j\sigma}(y')=(\mathbf B_{i}^2v)(y(y')),\quad
y'\in\gamma_{j\sigma}^\varepsilon.
$$

\begin{condition}\label{condB2vB1CConsistency}
For any function $v\in W_2^{2}(G\setminus\overline{\mathcal
O_{\varkappa_1}(\mathcal K)})$ and for any constant vector
$C=(C_1,\dots,C_N)$, the functions $B^v_{j\sigma}$ and $\mathbf
B_{j\sigma}C$, respectively, satisfy the consistency
condition~\eqref{eqConsistencyZ}.
\end{condition}

\begin{theorem}\label{thSmoothfne0}
Let Condition~$\ref{condProperEigen}$ be fulfilled. Then{\rm:}

1. If Condition~$\ref{condB2vB1CConsistency}$ holds and $u\in
W_2^1(G)$ is a generalized solution of problem~\eqref{eqPinG},
\eqref{eqBinG} with right-hand side $\{f_0,f_i\}\in L_2(G)\times
\mathcal S_2^{3/2}(\partial G)$, then $u\in W_2^2(G)$.

2. If Condition~$\ref{condB2vB1CConsistency}$ fails, then there
exists a generalized solution $u\in W_2^1(G)$ of
problem~\eqref{eqPinG}, \eqref{eqBinG} with certain right-hand
side $\{f_0,f_i\}\in L_2(G)\times \mathcal S_2^{3/2}(\partial G)$
such that $u\notin W_2^2(G)$.
\end{theorem}

We now consider problem~\eqref{eqPinG}, \eqref{eqBinG} with
homogeneous nonlocal conditions.

\begin{definition}\label{defAdmit}
We say that a function $v\in W_2^{2}(G\setminus\overline{\mathcal
O_{\varkappa_1}(\mathcal K)})$ is {\em admissible} if there exists
a constant vector $C=(C_1,\dots,C_N)$ such that
\begin{equation}\label{eqvadmissible}
B_{j\sigma}^v(0)+(\mathbf B_{j\sigma}C)(0)=0,\quad j=1,\dots,N,\
\sigma=1,2.
\end{equation}
Any vector $C$ satisfying relations~\eqref{eqvadmissible} is
called an {\em admissible vector corresponding to the function
$v$.}
\end{definition}

The set of admissible functions is linear. Clearly, the function
$v=0$ is an admissible function, whereas the vector $C=0$ is an
admissible vector corresponding to $v=0$. Moreover, one can verify
that any generalized solution of problem~\eqref{eqPinG},
\eqref{eqBinG} with homogeneous nonlocal conditions is an
admissible function.

Consider the following condition (which is weaker than
Condition~\ref{condB2vB1CConsistency}).

\theoremstyle{plain}
\newtheorem*{condition'}{Condition~\ref{condB2vB1CConsistency}$'$}
\begin{condition'}\label{condBv+BCConsist}
For any admissible function $v$ and for any admissible vector $C$
corresponding to $v$, the functions $B_{j\sigma}^v+\mathbf
B_{j\sigma}C$ satisfy the consistency
condition~\eqref{eqConsistencyZ}.
\end{condition'}

\newtheorem*{theorem'}{Theorem~\ref{thSmoothfne0}$'$}
\begin{theorem'}\label{thSmoothf0}
Let Condition~$\ref{condProperEigen}$ be fulfilled. Then{\rm:}

1. If Condition~$4'$ holds and $u\in W_2^1(G)$ is a generalized
solution of problem~\eqref{eqPinG}, \eqref{eqBinG} with right-hand
side $\{f_0,0\}$, $f_0\in L_2(G)$, then $u\in W_2^2(G)$.

2. If Condition~$4'$ fails, then there exists a generalized
solution $u\in W_2^1(G)$ of problem~\eqref{eqPinG}, \eqref{eqBinG}
with certain right-hand side $\{f_0,0\}$, $f_0\in L_2(G)$, such
that $u\notin W_2^2(G)$.
\end{theorem'}

The proofs of Theorems~\ref{thuinW_2^2NoEigen},
\ref{thSmoothfne0}, and~$2'$ are based on results concerning the
solvability of model nonlocal problems in plane angles in Sobolev
spaces~\cite{GurRJMP03} and on asymptotic behavior of solutions of
these problems in weighted spaces~\cite{SkMs86,GurPetr03}.

The author is grateful to Professor A. L. Skubachevskii for
attention to this work.

%
%

\end{document}